\documentclass[12pt,reqno]{amsart}	
\usepackage{amsthm,amsmath,amsfonts,amscd,amssymb,latexsym,exscale,xypic}
\usepackage[mathscr]{eucal}

\newtheorem{Theorem}{Theorem}[section]
\newtheorem{Lemma}[Theorem]{Lemma}
\newtheorem{Proposition}[Theorem]{Proposition}

\newtheorem{Scholium}[Theorem]{Scholium}

\newtheorem{Problem}[Theorem]{Problem}
\renewcommand{\qed}{\hspace*{\fill} $\Box$}
\newcommand{\tact}{| \! \bigcirc} 
\newcommand{\acc}{'\!}
\newcommand{\brac}[3]{(#1\,#2\,#3)} 

\newcommand{\demo}{{\sc Proof. \;\;}}
\newcommand{\complex}{{\mathbf C}} 
\renewcommand{\P}{{\mathbb P}}     
\renewcommand{\O}{{\mathcal O}}    
\newcommand{\G}{{\mathcal G}}      
\newcommand{\F}{{\mathcal F}}  
\newcommand{\Q}{{\mathcal Q}}  
\newcommand{\I}{{\mathcal I}}  
\newcommand{\J}{{\mathcal J}}  
\renewcommand{\aa}{{\mathcal A}}
\newcommand{\bb}{{\mathcal B}}

\newcommand{\hh}{{\mathcal H}}
\newcommand{\Hy}{{\mathbb H}}  
\newcommand{\TT}{{\mathbb T}}
\newcommand{\xbar}{\underline x} 
\newcommand{\ubar}{\underline u} 
\newcommand{\abar}{\underline a} 
\newcommand{\ybar}{\underline y} 
\newcommand{\vbar}{\underline v} 
\renewcommand{\emptyset}{\varnothing} 
\newcommand{\Sym}{\text{Sym}}
\renewcommand{\dim}{\text{dim}}
\newcommand{\codim}{\text{codim}}
\renewcommand{\ker}{\text{ker}}
\newcommand{\coker}{\text{coker}}
\newcommand{\rank}{\text{rank}}

\newcommand{\Res}{\text{Res}}
\newcommand{\image}{\text{image\,}}
\newcommand{\ra}{\rightarrow}
\newcommand{\RA}{\Rightarrow}
\newcommand{\lra}{\longrightarrow}

\newcommand{\incl}{\rto|<<\tip}

\begin{document}
\title{Decomposable Ternary Cubics} 
\author{Jaydeep V. Chipalkatti}
\maketitle 

\medskip 
\parbox{12cm}{\small
Cubic forms in three variables are parametrised by points 
of a projective space $\P^9$. We study the 
subvarieties in this space defined by decomposable forms. 
Specifically, we calculate their equivariant minimal resolutions 
and describe their ideals invariant-theoretically. }
\medskip \medskip 

Let $F$ be a nonzero cubic form in variables $x_1,x_2,x_3,$ written as 
follows: 
\begin{equation}
F = a_0\,x_1^3 + a_1\,x_1^2 x_2 + \dots + a_9\,x_3^3, 
\quad a_0,\dots,a_9 \in \complex. 
\label{F.expression} \end{equation}

We identify $F$ with the point $[a_0,\dots,a_9] \in \P^9$ and denote the homogeneous 
coordinate ring of $\P^9$ with $\complex[a_0,\dots,a_9]$. Now consider the 
set of those $F$ which factor as a product of a 
quadratic and a linear form, i.e., let 
\[ X_{\emptyset} = 
\{F \in \P^9: F = Q.L \;\; \text{for some forms $Q,L$ of degrees $2,1$}
\}. \] 
This is a projective subvariety of $\P^9$. The group 
$SL_3(\complex)$ acts on $\P^9$ as follows: 
given a matrix $A \in SL_3$, introduce new variables 
$[x_1',x_2',x_3'] = [x_1,x_2,x_3]A$, and 
define $a_0',\dots,a_9'$ by forcing the identity
\[ a_0\,x_1^3 + a_1\,x_1^2 x_2 + \dots + a_9\,x_3^3 = 
   a_0'\,{x_1'}^3 + a_1'\,{x_1'}^2 x_2' + \dots + a_9'\,{x_3'}^3.  \] 
Then $[a_0,\dots,a_9] \stackrel{A}{\lra} [a_0',\dots,a_9']$. The 
imbedding $X_{\emptyset} \subseteq \P^9$ is thus $SL_3(\complex)$--equivariant. 
Hence $I_{X_\emptyset} < \complex[a_0,\dots,a_9]$ is a 
representation of $SL_3$, and so are all the syzygy modules 
in the minimal resolution of $I_{X_\emptyset}$. 

There are different possible factorisations of $F$, leading to five similarly defined
varieties:
\[ \begin{array}{ccl} 
X_{\equiv} & = & 
\{F \in \P^9 : F = L^3 \;\; \text{for a linear form $L$} \}, \\ 
X_{\neq} & = & 
\{F \in \P^9 : F = L_1^2L_2 \;\; \text{for some $L_i$} \}, \\ 
X_{\Delta} & = & 
\{F \in \P^9 : F = L_1L_2L_3 \;\; \text{for some  $L_i$} \}, \\ 
\end{array} \]
and the loci $X_{\tact},X_Y$ defined as Zariski closures
\begin{equation} \begin{array}{ccl} 
X_{\tact} & = & 
\{F \in \P^9 : F = Q.L \;\; \text{where the conic $Q=0$ is smooth and } \\
 & & \; \quad \text{ tangent to the line $L=0$} \overline{\} \,}, \\
X_Y & = & 
\{F \in \P^9 : F = L_1L_2L_3 \;\; \text{where $L_i=0$ are 
concurrent lines} \overline{\} \, }. 
\end{array} \end{equation}
All the loci are irreducible, and there are inclusions 
\[ \diagram
X_{\equiv}^{(2)} \incl & X_{\neq}^{(4)} \incl & 
X_{Y}^{(5)} \incl \dto|<<\tip & X_{\Delta}^{(6)} \dto|<<\tip   \\ 
 & & X_{\tact}^{(6)} \incl & X_{\emptyset}^{(7)} 
\enddiagram \] 
The superscripts indicate dimensions. 
Consider the following two problems: 
\begin{itemize} 
\item 
[(I)] calculating the $SL_3$-equivariant minimal resolution for each variety $X$, and 
\item 
[(II)] expressing the minimal set of generators of $I_X$ (i.e., the degree 
zero syzygies), in the language of classical invariant theory. 
\end{itemize} 

Problem II, in a slightly weaker form, is of classical origin. A priori, we 
are looking for a finite set of concomitants of $F$, such that 
their vanishing is a necessary and sufficient condition for 
$F$ to lie in $X$. Describing such a set is tantamount to describing $I_X$ upto 
its radical. 

This is an instance of the general 
philosophy that for a form $F$ (of any degree in any number 
of variables), any property which is invariant under a linear change of 
variables (e.g., decomposabililty, existence of singular points, expressibility 
as a sum of a specified number of powers of linear forms) should be 
characterisable by the vanishing (or non-vanishing) of certain 
concomitants of $F$. Classical literature on invariant theory contains a 
very large number of results along these lines--see 
e.g.~\cite{DolgachevKanev,GrYo,Salmon2,Weitzenbock} or the collected 
papers of Cayley and Clebsch. Salmon's book \cite{Salmon1} is an excellent 
reference for the classical theory of plane curves and connections to invariant 
theory. Its chapter V is entirely devoted to ternary cubics. For newer (i.e.
post-Hilbert) directions in invariant theory, see \cite{CD,Kraft,Weyl}. 

Problem I is a natural generalisation of II, but certainly more modern. In this 
paper we solve problem I, and then interpret the representations corresponding to 
the ideal generators as concomitants of ternary cubics.  

Some of our varieties have alternate descriptions, e.g., $X_\equiv$ is 
the cubic Veronese imbedding of $\P^2$ and $X_Y$ is a rank variety in the sense of 
Porras \cite{Porras}. As such, the minimal systems of ideal generators for them 
are known. The standard equations defining the Veronese are given in 
\cite[p.~23]{Harris}, but not in invariant-theoretic language. 
We will describe Porras's deduction of the minimal resolution of 
$X_Y$ in \S \ref{section.XY}. In 
the book referred to above, Salmon deduces sets of equations describing 
the varieties $X_\Delta$ and $X_\emptyset$ set-theoretically.
Our calculations confirm that these equations generate the ideal 
$I_X$ in each case. 

The next section contains background material on representation theory. For each of the 
loci $X$, we find the ideal $I_X$ and its Betti numbers using 
a machine computation. Then we identify the syzygy modules as $SL_3$-representations 
using the hypercohomology spectral sequence, together with algebraic
considerations specific to each $X$. In the last section, we write down the 
concomitants corresponding to the ideal generators. 

This paper was inspired by the \emph{Weltanschauung} of 
\cite{FH} (especially \S 11.3, \S 13.4), and also by Salmon's 
wonderful book \cite{Salmon1}. The program Macaulay-2 has been 
of immense help in the computations, and it is a pleasure
to thank its authors Daniel Grayson and Michael Stillman. 
I thank the referee for helpful observations and Anthony V. Geramita, 
Leslie Roberts and the Queen's University for financial assistance 
during the course of this work. 

\setcounter{equation}{0}

\section{Preliminaries} 
The purpose of this section is to establish notation and state the working 
definitions, all of which are phrased in invariant-theoretic 
terms. I have written somewhat expansively, in the hope that the 
framework introduced here may have more general application. 

Let $V$ be a three dimensional $\complex$--vector space. All 
representations considered henceforth are for the group $SL(V)$. 
(See \cite{Fulton,FH,Sturmfels} for the theory.) 
For a pair $(m,n)$ of nonnegative integers, we have the Schur 
module $S_{m+n,n} = S_{m+n,n}(V)$, with the convention that 
$S_m = S_{m,0} = \Sym^m(V)$ and $S_{1,1} = \wedge^2 V$. Note the 
isomorphism\footnote{Henceforth `$=$' means an isomorphism of 
$SL$-modules.} $S_{m+n,n}^* = S_{m+n,m}$ and the formula 
\begin{equation}
\dim \, S_{m+n,n} = \frac{1}{2}(m+1)(n+1)(m+n+2). 
\end{equation}

In the sequel, we frequently need to decompose the tensor product 
$S_{m_1+n_1,n_1} \otimes S_{m_2+n_2,n_2}$ and 
the plethysm $S_\lambda(S_{m+n,n})$ as direct sums
of irreducible representations. The former is governed by 
the Littlewood-Richardson rule and the latter 
by formulae (A.5), (A.21) in \cite[appendix A]{FH}. This was programmed in 
Maple by the author.\footnote{The `Symmetrica' package and 
John Stembridge's `SF' package do this too, but I have found 
that ad-hoc routines for dim $V=3$ tend to work 
faster. Of course the issue of speed is more serious for the 
plethysm.}

Let $\xbar  = \{x_1,x_2,x_3\}$ be a basis of $V^*$. The 
identification $\wedge^2 V^* = V$ gives a basis 
$\ubar = \{u_1,u_2,u_3\}$ of $V$, 
where 
\begin{equation}
u_1 = x_2 \wedge x_3,\; u_2 = x_3 \wedge x_1,\; u_3 = x_1 \wedge x_2. 
\end{equation}
If $a_r$ appears as the coefficient of $x_1^{i_1}x_2^{i_2}x_3^{i_3}$ in 
(\ref{F.expression}), redefine $a_r$ to be 
$u_1^{i_1}u_2^{i_2}u_3^{i_3}$. Thus $ \abar = \{ a_0,\dots,a_9 \}$ is a basis 
of $S_3$ and the form $F$ is recast as the canonical trace element 
\[ \sum\limits_{i_1+i_2+i_3=3}
 u_1^{i_1}u_2^{i_2}u_3^{i_3} \otimes x_1^{i_1}x_2^{i_2}x_3^{i_3}
 \quad \in S_3 \otimes S_3^*. \] 
The ring $\complex[a_0,\dots,a_9]$ is identified with 
the symmetric algebra $R=\bigoplus\limits_{\ell \ge 0} S_\ell\, ( S_3)$, 
and we declare $\P S_3^* = \text{Proj}\, R$ to be the ambient 
space for all loci under consideration. 

\subsection{A Basis for $S_{m+n,n}$} \label{basis.Smn}
There is a canonical injection (see [p.233 ff, loc.cit.])
\[ S_{m+n,n} \lra \Sym^n(\wedge^2 V) \otimes \Sym^m(V), \] 
and the right hand side has a basis of monomials 
\begin{equation}
\{x_1^{n_1}x_2^{n_2}x_3^{n_3} \otimes 
u_1^{m_1}u_2^{m_2}u_3^{m_3}: \;\Sigma n_i = n,\, \Sigma m_i = m \}.  
\label{monomiallist} \end{equation}
Let $T$ be a semistandard tableau on numbers $1,2,3$ on the Young 
diagram of $(m+n,n)$. Say the entries in the first (resp. second) row are 
$ r_1 \le \dots \le r_n \le s_1 \le \dots \le s_m$ (resp. 
$t_1 \le \dots \le t_n$), with $r_i < t_i$ for $1 \le i \le n$. 
Let $x_{\langle r_i,t_i \rangle}$ denote resp.~$x_1,-x_2,x_3$ 
if $(r_i,t_i)=(2,3),(1,3),(1,2)$, and write 
\begin{equation}
 X_T = \prod\limits_{i=1}^n x_{\langle r_i,t_i \rangle} 
 \otimes \prod\limits_{i=1}^m u_{s_i}
\label{basislist} \end{equation}
(E.g., for the tableau 
$T = [1\,2\,2\,2\,3|2\,3\,3]$ on $(5,3)$, 
$X_T = x_1^2x_3 \otimes u_2u_3$.) 
Now (ignoring signs) the subset 
$\{X_T : T \; \text{semistandard} \}$ of 
(\ref{monomiallist}) is a basis for $S_{m+n,n}$. 
\subsection{Concomitants} \label{section.concomitants}
Let us assume that we have 
an injective map of representations 
\[ S_{m+n,n} \stackrel{\varphi}{\lra} S_\ell\,(S_3), \] 
or equivalently, an imbedding 
\[ \complex \lra S_\ell\,(S_3) \otimes S_{m+n,n}^* 
(= S_\ell\,(S_3) \otimes S_{m+n,m}). \] 
Let $\Phi$ denote the image of \,$1$, which is a ${\mathbf Q}$--linear 
combination of monomials 
\[ \{\, \prod\limits_{i=0}^9 a_i^{l_i} \otimes 
     \prod\limits_{i=1}^3 x_i^{m_i} \otimes 
     \prod\limits_{i=1}^3 u_i^{n_i} : \; 
\Sigma\, l_i = \ell,\, \Sigma\, m_i = m,\, \Sigma\, n_i = n \}.  \] 
In 19th century terminology (due to Sylvester), $\Phi$ is called 
a concomitant of degree $\ell$, order $m$ and class $n$ (for 
ternary cubics). We indicate this by writing $\Phi_{\ell,m,n}$. 
Of course, there may exist more than one concomitant for 
given $(\ell,m,n)$ (or none at all). 
A concomitant is called a covariant (resp. contravariant)
if its class (resp. order) is zero, and an invariant if 
both are zero. Note that the imbedding $\varphi$ can be 
recovered from $\Phi$. 

For instance, the Hessian of a ternary cubic is a covariant 
of degree $3$ and order $3$. The dual of the curve 
$F=0$ is defined by a contravariant of degree $4$ and class $6$. 
\subsection{The Symbolic Method}
We will represent concomitants using the German symbolic method. 
A brief explanation follows, but it is unlikely to 
be intelligible without prior acquaintance with the method. 
The notation follows Grace and Young \cite{GrYo}, also 
see \cite[Ch.~6]{Olver} for a modern exposition. 

We will use the letters $x_1,x_2,x_3; u_1,u_2,u_3$. In addition, 
we have an indefinite supply of indexed Greek letters 
$\alpha_1,\alpha_2,\alpha_3; \beta_1,\beta_2,\beta_3$ etc. Set 
\[ 
\alpha_x = \alpha_1 x_1 + \alpha_2 x_2 + \alpha_3 x_3 \quad 
\text{and ditto for $\beta_x,\gamma_x$ etc}. \] 
Formally write 
$ F = \alpha_x^3 = \beta_x^3 = \dots $. That is to say, after 
expansion, 
$\frac{3!}{i_1!i_2!i_3!}\alpha_1^{i_1}\alpha_2^{i_2}\alpha_3^{i_3}$ 
stands for the appropriate coefficient $a_r$ in 
(\ref{F.expression}). 
The symbols $\brac{\alpha}{\beta}{\gamma}$, $\brac{\alpha}{\beta}{u}$ 
respectively stand for the determinants 
\[ \left| \begin{array}{ccc} 
\alpha_1 & \alpha_2 & \alpha_3 \\ 
\beta_1 & \beta_2 & \beta_3 \\ 
\gamma_1 & \gamma_2 & \gamma_3 \end{array} \right| \quad \text{and}\quad 
\left| \begin{array}{ccc} 
\alpha_1 & \alpha_2 & \alpha_3 \\ 
\beta_1 & \beta_2 & \beta_3 \\ 
u_1 & u_2 & u_3 \end{array} \right|. \] 

Now for instance, the Hessian is written 
$\brac{\alpha}{\beta}{\gamma}^2\alpha_x \beta_x \gamma_x$. To 
evaluate this, multiply out the expression formally, and substitute 
the $a_r$'s as above. The result, appropriately, is a polynomial 
of degree $3$ each in the $\abar,\xbar$. As another example, the 
dual curve mentioned above has equation 
${\brac{\alpha}{\beta}{u}}^2 {\brac{\gamma}{\delta}{u}}^2
\brac{\alpha}{\gamma}{u}\brac{\beta}{\delta}{u}=0$.

Each concomitant $\Phi_{\ell,m,n}$ has a (non-unique) 
expression as a sum of products of symbols of the form 
$\alpha_x, \brac{\alpha}{\beta}{\gamma}, \brac{\alpha}{\beta}{u}$. 
In \emph{each} summand, the letters $x,u$ should occur resp.~$m,n$ times, 
along with $\ell$ Greek letters, each occuring exactly thrice. 
Conversely, every such symbolic expression corresponds to a 
concomitant, unless it vanishes identically after substitution. 

\subsection{A Spectral Sequence} 
Let $X$ be any of the loci above. The equivariant minimal 
resolution of its defining ideal $I_X < R$, will be written 
\[  \dots \ra F^{\,p} \ra F^{\,p+1} 
\ra \dots \ra F^{\,0} \ra I_X \ra 0, \] 
\begin{equation} \text{with}\quad 
 F^{\,p} = \bigoplus_{j > 0 } M_{j,-p} \otimes R(-j+p) 
\quad \text{for $ p \le 0$.} \label{ideal.RES} \end{equation}

Letting $\F^{\,p} = (F^{\,p}\widetilde{)} = 
\bigoplus M_{j,-p} \otimes \O_{\P^9}(-j+p)$, we have a resolution 
$\F^{\,\bullet} \ra \I_X \ra 0$, for the ideal sheaf. For 
$\ell \in {\mathbf Z}$, let $\F^{\,p}(\ell) =\F^{\,p} \otimes \O_{\P}(\ell)$. 
Then we have a second quadrant spectral sequence 
\begin{equation}
\begin{array}{ll}
{\check E}_1^{p,q} = H^q(\P^9,\F^{\,p}(\ell)) 
& \text{for $p \le 0,\, q \ge 0$,}  \\
d_r: {\check E}_r^{p,q} \ra {\check E}_r^{p+r,q-r+1}; & 
\;{\check E}_\infty^{p,q} \RA H^{p+q}(\P^9,\I_X(\ell)).
\end{array} \label{ideal.spectralseq} \end{equation}
This will be used in conjunction with the standard calculation 
of cohomology of line bundles on projective space 
(\cite[Ch. III, Theorem 5.1]{Ha}) and Serre duality. 

The bulk of the paper is concerned with the identification of 
the $M_{j,-p}$ \emph{qua} $SL_3$--representations. The ideals 
and their Betti numbers (i.e. dimensions of $M_{j,-p}$) 
were calculated in Macaulay-2. 

\subsection{Computation of $I_X$.}
For illustration, consider the locus $X_{\neq}$. It is the image of 
the projective morphism 
\[ \P V^* \times \P V^* \lra \P S_3^* , \quad (L_1,L_2) \ra L_1^2L_2. \]
Let 
$ L_1 = b_1x_1+b_2x_2+b_3x_3, \quad L_2 = c_1x_1+c_2x_2+c_3x_3$,
where the $b_i,c_i$ are indeterminates. 
By forcing the equality $F = L_1^2L_2$, we get polynomial expressions 
\[ 
a_r = \varphi_r(b_1,b_2,b_3; c_1,c_2,c_3) 
\quad \text{for $r=0,\dots,9$};  \]
defining a ring map 
$f_{\neq}: \complex[a_0,\dots,a_9] \lra \complex[b_1,\dots,c_3]$. 
Then $I_{X_{\neq}}$ equals $\ker (f_{\neq})$. Its minimal 
resolution is tabulated as follows: 
\[
\begin{tabular}{c|ccccccc} 
  & 0  & 1   &  2  & 3   & 4  & 5  & 6  \\ \hline 
3 & 20 & 45  & 36  & 10                 \\ 
4 & 28 & 126 & 225 & 200 & 90 & 18 & 1 
\end{tabular}  \]
In the notation of (\ref{ideal.RES}), 
the entry in row $j$ and column $-p$ is dim $M_{j,-p}$. 
E.g., $M_{3,2}$ is $36$--dimensional. 

The loci $X_{\equiv},X_{\Delta},X_{\emptyset}$ are images of 
projective morphisms, respectively from $\P V^*, (\P V^*)^3, 
\P S_2^* \times \P V^*$ to $\P S_3^*$; hence a similar procedure 
works for them. That these loci are closed, follows from the fundamental 
theorem of elimination theory. 

 In order to find $I_{X_{\tact}}$, we need the tact-invariant of 
a conic and a line (see \cite[\S 96]{Salmon1}). Let 
\[ Q = q_1x_1^2 + q_2x_2^2 +q_3x_1x_2+q_4x_1+q_5x_2+q_6, \;\;
   L = b_1x_1 + b_2x_2 +b_3. \]
Writing $Q.L = F|_{(x_3=1)}$, we get expressions 
$a_r = \varphi_r(q_1,\dots,b_3)$ as before. 
Let $\Res = \text{Resultant}(Q,L;x_2)$, then the $x_1$-values 
for which $\Res=0$ are the $x_1$-coordinates of the points 
$Q=L=0$. The condition that they coincide is 
$T' =\text{Discriminant}(\Res,x_1)=0$. Then the tact-invariant 
$T = \frac{1}{b_2^2} T'$. It vanishes iff the line $L=0$ is tangent to $Q=0$. 
(The extraneous factor appears  because when $b_2=0$, the 
$x_1$-coordinates coincide for any position of $Q$.) For what it is worth, 
\[ \begin{aligned} T  = & \, 
 4q_2 q_6 b_1^2 - 4q_2q_4b_1b_3 + 4q_1q_2b_3^2 -
  q_5^2b_1^2 - 4q_3q_6b_1b_2 + 2q_4q_5b_1b_2 + \\ 
  & \,  2q_3q_5b_1b_3 - q_4^2b_2^2 - 4q_1q_5b_2b_3 + 2q_3q_4b_2b_3
  - q_3^2b_3^2 + 4q_1q_6b_2^2. 
\end{aligned} \] 

Now $I_{X_{\tact}}$ is the kernel of the map 
\[ f_{\tact}: 
\complex[a_0,\dots,a_9] \lra \complex[q_1,\dots,b_3]/(T). \]

To calculate $I_{X_Y}$, let 
\[ \begin{array}{c} 
   L_1 = b_1x_1 + b_2x_2 +b_3x_3 \\ 
   L_2 = c_1x_1 + c_2x_2 +c_3x_3  \\
   L_3 = d_1x_1 + d_2x_2 +d_3x_3 \\ 
\end{array} \; \text{and} \;\; D = 
\left| \begin{array}{ccc} 
 b_1 & b_2 & b_3 \\ 
 c_1 & c_2 & c_3 \\ 
 d_1 & d_2 & d_3 \end{array} \right|. \] 
Then $D=0$ is the condition that the lines be concurrent. 
The $\varphi_r$ are given by $F=L_1L_2L_3$, and $I_{X_Y}$ is 
the kernel of 
\[ f_Y: \complex[a_0,\dots,a_9] \lra \complex[b_1,\dots,d_3]/(D). \]

\subsection{Characters} 
We write the formal character of $S_{a,b}$ as $[S_{a,b}] = \sigma_{a,b}$. 
Every finite dimensional $SL_3$-module $U$ is uniquely expressible as 
a direct sum $\bigoplus (S_{a_i,b_i})^{\oplus n_i}$, 
hence write $[U] = \sum n_i\, \sigma_{a_i,b_i}$. 

If $C^\bullet = \{C^p\}_p$ is a bounded complex of finite dimensional 
$SL_3$--representations and $SL$-equivariant differentials, then 
define 
$[C^\bullet] = \sum\limits_p (-1)^p\, [C^p]$. 
By Schur's lemma, we have 
\begin{equation}
 [C^\bullet] = \sum\limits_p (-1)^p\, [H^p(C^\bullet)].
\label{complex.char} \end{equation}
\setcounter{equation}{0}

\section{The Locus $X_{\equiv}$}
The Betti numbers of $I_X$ are 
\begin{equation} \begin{tabular}{c|ccccccc}
  &  0 & 1   &   2 &   3 &   4 &  5 & 6 \\ \hline 
2 & 27 & 105 & 189 & 189 & 105 & 27 &    \\
3 &    &     &     &     &     &    & 1 
\end{tabular} \label{betti.equiv} \end{equation}

Now $X$ is the cubic Veronese imbedding of $\P^2$, hence 
arithmetically Cohen--Macaulay. Since $\omega_X = \O_X(-1)$, it is 
moreover subcanonical, hence arithmetically Gorenstein. Thus the 
minimal resolution of $R/I_X$ is centrally symmetric, i.e., 
$M_{23},M_{24},M_{25} \footnotemark$ are 
\footnotetext{\,Henceforth the commas are left out.}
respectively dual to $M_{22},M_{21},M_{20}$.

Consider the diagram 
\[ \begin{array}{cll}
0 \lra  H^0(\I_X(\ell)) \lra & H^0(\P^9,\O_{\P}(\ell)) & 
\stackrel{\gamma}{\lra} H^0(X,\O_X(\ell)) \\
 & \quad \;\;  \big| \big| & \qquad \quad \big| \big| \\
 & S_\ell(S_3) & \lra  \;\; S_{3\ell}
\end{array} \]
The map $\gamma$ is surjective for $\ell \ge 0$, 
hence by (\ref{complex.char}), 
\begin{equation} [H^0(\I_X(\ell))] = [H^0(\O_{\P}(\ell))] - 
[H^0(\O_X(\ell))]. 
\label{E1.equiv} \end{equation}
Substitute $\ell=2,3$, and decompose the plethysms 
involved, then
\[ [H^0(\I_X(2))] = \sigma_{42}, \quad 
   [H^0(\I_X(3))] = \sigma_{72}+ \sigma_{63}+ \sigma_{33}+ \sigma_{30}. 
\footnotemark \] 
\footnotetext{\dots and even here.} 
Now successively let $\ell=2,3$ in the spectral sequence 
(\ref{ideal.spectralseq}). In each case, all nonzero 
terms are in the row $q=0$, so 
${\check E}_2 = {\check E}_\infty$. Thus we have 
$M_{20} = H^0(\I_X(2))$ and an exact sequence 
\[ 0 \lra M_{21} \lra M_{20} \otimes S_3 \lra H^0(\I_X(3)) \lra 0. \] 
It follows that $M_{20} = S_{42}$ and  
$M_{21} = S_{54}\oplus S_{51} \oplus S_{42} \oplus S_{21}$; 
hereafter written $M_{21} = \{54,51,42,21\}$. Now 
$M_{22}$ is calculated similarly using $\ell=4$. This 
is the end result: 
\begin{equation} 
 M_{20}  = \{ 42\}, M_{21} = \{ 54,51,42,21 \}, 
 M_{22} = \{ 63,54,51,42,33,30,21\}, 
\label{syzygymodules.equiv} \end{equation}
and since these are all self-dual, $M_{25} = M_{20},
M_{24} = M_{21},M_{23} = M_{22}$ and of course $M_{36}=\{00\}$. 

\begin{Problem} \rm 
All the syzygy modules in the resolution are self-dual. 
One would like to know if this has any geometric significance and to what 
extent this is true of other Veronese imbeddings. (This is trivially 
true of the rational normal curve, since all $SL_2$-representations are 
self-dual. But it fails for the quadratic imbedding of $\P^2$ in $\P^5$.) 
\end{Problem} 

As far as I know, the problem of describing an equivariant minimal 
resolution for the Veronese variety is open in general. See \cite{Green} 
for some results. 

\setcounter{equation}{0}

\section{The Locus $X_{\neq}$} 
The Betti numbers of $I_X$ are
\begin{equation}
\begin{tabular}{c|ccccccc}
  & 0  & 1   &  2  & 3   & 4  & 5  & 6  \\ \hline
3 & 20 & 45  & 36  & 10                 \\ 
4 & 28 & 126 & 225 & 200 & 90 & 18 & 1 
\end{tabular} \label{betti.neq} \end{equation} 
We have an exact sequence 
\begin{equation} 
0 \lra \I_{X_{\neq}} \lra \O_{\P^9} \lra \O_{X_{\neq}} \lra 0, 
\label{defn.IX} \end{equation}
and to use the spectral sequence (\ref{ideal.spectralseq}), some knowledge 
of the groups $H^\bullet(\P^9, \I_X(\ell))$ is needed. We 
proceed in several steps. 
\subsection{Step I} 
Consider the commutative triangle 
\[ \diagram 
\P V^* \dto_{\delta} \drto^g \\ 
\P V^* \times \P V^* \rto_(0.6){f}  & \P S_3^* 
\enddiagram \] 
where $f: (L_1,L_2) \ra L_1^2L_2$, and $\delta$ is the diagonal 
imbedding. Then $\image f = X_{\neq}, \image g = X_{\equiv}$. 
Define an $\O_{\P^9}$-module $\Q$ as the cokernel 
\begin{equation}
0 \lra \O_{X_{\neq}} \lra f_* \O_{\P^2 \times \P^2}
\lra \Q \lra 0, \label{defn.Q} \end{equation} 
then $\text{supp}(\Q) = X_{\equiv}$.

\begin{Lemma} There is an isomorphism 
$g^* \Q = \Omega^1_{\P^2}$. 
\label{lemma.calc.Q} \end{Lemma}
\demo Let $\J$ be the ideal sheaf of image$(\delta)$, thus 
\[ 0 \lra \J \lra \O_{\P^2 \times \P^2} \ra 
\delta_* \O_{\P^2} \lra 0. \]
We have a commutative ladder 
\[ \diagram
0 \rto & \ker\, 1 \dto^2 \rto & \O_{X_{\neq}} \dto^3 \rto^1 & 
\O_{X_{\equiv}} \dto^4 \rto & 0 \\ 
0 \rto & f_*(\J) \rto & f_* \O_{\P^2 \times \P^2} \rto & 
g_* \O_{\P^2} \rto & 0 
\enddiagram \] 
(Since $f$ is a finite morphism, $f_*$ is exact.) 
Since $4$ is an isomorphism, $\coker\, 2 = \coker\, 3$, 
giving a surjection 
$f_*(\J) \stackrel{q}{\twoheadrightarrow} \Q$. By 
(\cite[p.110]{Ha}), we have an adjoint morphism 
$\J \stackrel{q'}{\lra} f^* \Q$. Now the composite 
\[ f^*f_*(\J) \lra \J \stackrel{q'}{\lra} f^* \Q \] 
is surjective, since $f^*$ is right-exact. Hence $q'$ must be 
surjective. Applying $\delta^*$, we get a surjection 
$\delta^*(\J) \stackrel{q''}{\twoheadrightarrow} g^* \Q$. 
By definition, $\Omega^1_{\P^2}$ equals $\delta^*(\J)$. 

To show that $q''$ is an isomorphism, it suffices to 
observe that both sheaves have the same Hilbert polynomial. 
By the K{\"u}nneth formula, 
\[ H^0(f_* \O_{\P^2 \times \P^2} \otimes \O_{\P^9}(\ell)) = 
H^0(\P^2,\O_{\P}(2 \ell)) \otimes H^0(\P^2,\O_{\P}(\ell)),
\] and this combined with 
(\ref{betti.neq}),(\ref{defn.Q})  gives 
$9\ell^2 -1$ as the Hilbert polynomial of $\Q$. As 
$g^* \O_{\P^9}(\ell) = \O_{\P^2}(3 \ell)$, that of 
$g^* \Q$ is $\ell^2 -1$. That this is also the Hilbert 
polynomial of $\Omega^1_{\P^2}$, can be seen from the Euler sequence. 
The lemma is proved. \qed 

\subsection{Step II}
The next step is to construct a complex of $SL_3$-modules 
with \emph{explicitly computable terms}, such that the global 
sections of $\I_X(\ell)$ and $\Q(\ell)$ appear as its cohomology 
modules. We will use the Borel--Weil--Bott theorem 
several times (see e.g. \cite[Theorem 0.1]{ego1}). 

Let $Y$ be the variety $\P V^* \times \P V^* \times \P S_3^*$ with 
projections $\mu_1,\mu_2,\pi$ onto the successive factors. 
Consider the diagram 
\[ \diagram 
\P V^* \times \P V^* \rto^(0.65){i} \drto_f & Y \dto^{\pi} & 
\quad (L_1,L_2) \rto^(0.45){i} & (L_1,L_2,L_1^2L_2) \\ 
 & \P S_3^* \enddiagram \] 
and let $\Gamma \subseteq Y$ be the image of $i$. Then 
$\pi(\Gamma) = X_{\neq}$ (henceforth written $X$). We will 
derive a Koszul resolution of $\O_\Gamma$. 
Define a  vector bundle 
\[ \aa = \mu_1^* \O_{\P^2}(2) \otimes \mu_2^* \O_{\P^2}(1)
\otimes \pi^* T_{\P^9}(-1) \]
on $Y$. If $y = (L_1,L_2,F)$ be a point of $Y$, then the fibre 
of $\aa$ over $y$ is the $9$-dimensional vector space 
\[ \langle L_1^2 \rangle^* \otimes \langle L_2 \rangle^* \otimes  
 {S_3^*} / {\langle F \rangle}. 
\] 
Here $\langle L_1^2 \rangle \subseteq S_2^*$ is the one dimensional 
space generated by $L_1^2$, etc. 

By the K{\"u}nneth formula,  
$H^0(Y,\aa) = S_2 \otimes S_1 \otimes S_3^*$. 
The multiplication map $S_2^* \otimes S_1^* \ra S_3^*$ corresponds 
to a distinguished section $\xi$ in $H^0(Y,\aa)$. 
\begin{Lemma}
The zero scheme of $\xi$ is $\Gamma$. 
\end{Lemma} 
\demo To say that $\xi$ vanishes at $y$, is to say that the 
element $L_1^2L_2 + \langle F \rangle 
\in {S_3^*} / {\langle F \rangle}$ is zero. This happens iff 
$y \in \Gamma$. \qed 

Now $\codim(\Gamma,Y) = \rank \, \aa = 9$, hence the Koszul 
complex of $\xi$ resolves $\Gamma$. 
\begin{equation} \begin{aligned} 
0 \ra \wedge^{9}\aa^* \ra \dots \wedge^{-p} \aa^*
\stackrel{\wedge \xi}{\lra}  \wedge^{-(p+1)} \aa^* \ra \dots \O_Y
& \ra  \O_\Gamma \ra 0, \\ 
& \text{for $ -9 \le p \le 0$.} \end{aligned} 
\end{equation} 
Consider the hypercohomology spectral sequence 
\begin{equation} \begin{array}{ll}
 E_1^{p,q}= {\mathbf R}^q\pi_*(\wedge^{-p} \aa^*), & \quad 
d_r: E_r^{p,q} \lra E_r^{p+r,q-r+1}; \\
 E_\infty^{p+q} \RA {\mathbf R}^{p+q}\pi_*(\O_\Gamma) 
& \quad \text{for $-9 \le p \le 0,\, q \ge 0$.}
\end{array} \label{neq.SP2} \end{equation}
Now 
\[ \wedge^{-p} \aa^* = \mu_1^* \O_{\P^2}(2p) \otimes \mu_2^* \O_{\P^2}(p)
\otimes \pi^* \Omega^{-p}_{\P^9}(-p); \]
hence 
\[ {\mathbf R}^q\pi_*(\wedge^{-p} \aa^*) =  
\{ \bigoplus_{i+j=q} H^i(\P^2, \O_{\P^2}(2p)) \otimes H^j(\P^2, \O_{\P^2}(p)) 
\} \otimes \Omega^{-p}_{\P^9}(-p). \]

Apart from $E_1^{0,0}=\O_{\P^9}$, all nonzero terms are 
concentrated in the row $q=4$ and columns $-9 \le p \le -3$. 
Using Serre duality,
\begin{equation} \begin{aligned} E_1^{p,4} & = 
H^2(\P^2,\O_{\P^2}(2p)) \otimes H^2(\P^2,\O_{\P^2}(p))
\otimes \Omega^{-p}(-p) \\
& = S_{-(2p+3)}^* \otimes  S_{-(p+3)}^* \otimes \Omega^{-p}(-p).
\label{neq.EXPQ} \end{aligned} \end{equation}
Define a complex $\G^\bullet$ by $\G^{\,p} = E_1^{\,p,4}$, with $d_1$ 
as the differential. It lives in the range 
$-9 \le p \le -3$. 

\begin{Proposition} The cohomology of $\G^\bullet$ equals 
\[ \hh^p(\G^\bullet) ( = E_2^{p,4} ) = 
\begin{cases} 
\I_X  & \text{if $p=-5$,} \\
\Q  & \text{if $p=-4$,} \\ 
\, 0 & \text{otherwise.} \end{cases} \]
\end{Proposition}

\demo Clearly $E_2 = \dots = E_5$ and $E_6 = E_\infty$. 
Since $\pi|_\Gamma$ is a finite morphism, 
${\mathbf R}^{\ge 1}\pi_* \O_{\Gamma} = 0$. Hence 
for $p \neq -5,-4$, we have $E_2^{p,4} = E_\infty^{p,4}=0$.
Consider the extension 
\[ 0 \ra E_6^{0,0} \ra \pi_* \O_\Gamma \ra E_6^{-4,4} \ra 0.\]
The term $E_6^{0,0}$ is the image of the natural map 
$E_1^{0,0} (= \O_{\P^9}) \ra \pi_* \O_\Gamma$, which is $\O_X$. 
Hence 
\[ \begin{aligned} 
E_5^{-5,4} & = \ker\, (\O_{\P^9} \ra \O_X) & = \I_X, \\ 
E_5^{-4,4} & = \coker\, (\O_X \ra \pi_* \O_\Gamma) & = \Q. 
\end{aligned} \]  \qed 
 
For $\ell \in {\mathbf Z}$, write $\G^\bullet(\ell) = 
\G^\bullet \otimes \O_{\P^9}(\ell)$. All further calculation 
revolves around the hypercohomology of $\G^\bullet(\ell)$. 
There are two second quadrant spectral sequences 

\begin{equation} \begin{alignat}{3}
& E_1^{p,q} = H^q(\P^9,\G^p(\ell)),\quad 
& d_r  :  E_r^{p,q}  & \lra E_r^{p+r,q-r+1};  \label{neq.SP3} \\
& \acc E_2^{p,q} = H^q(\P^9,\hh^p(\G^\bullet(\ell))), \quad 
& \acc d_r  : \; \acc E_r^{p,q}  & \lra\; \acc E_r^{p-r+1,q+r};
\label{neq.SP4} \end{alignat} \end{equation}
with 
$E_\infty^{p,q},\; \acc E_\infty^{p,q} \RA \Hy^{p+q}(\G^\bullet(\ell))$. 
If $\ell$ is clear from the context, we write $\Hy^j$ for 
$\Hy^j(\G^\bullet(\ell))$. 
Henceforth a symbol such as $E_r^{p,q}$ refers to (\ref{neq.SP3}), 
and \emph{not} to (\ref{neq.SP2}). 

\subsection{Step III: Analysis of (\ref{neq.SP3})} The terms 
\begin{equation}
 E_1^{p,q} = S^*_{-(2p+3)} \otimes S^*_{-(p+3)} \otimes 
H^q(\P^9, \Omega^{-p}_{\P^9} (\ell-p)) 
\label{E1pq.EXP} \end{equation}
and $\acc E_2^{-4,q} = H^q(\P^9,\Q(\ell))$ 
are evaluated by appealing to the Borel--Weil--Bott theorem. 
The result is as follows: 
\begin{enumerate} 
\item{Assume $\ell >0$, and let 
$\lambda =(\ell,\underbrace{1,\dots,1}_{\text{$-p$ times}})$. Then 
\[ H^q(\Omega^{-p}(\ell-p)) = 
\begin{cases} 
S_\lambda(S_3) & \text{for $q=0$,} \\ 
0            & \text{for $q \neq 0$.} \end{cases} \] 
Moreover $H^q(\Q(\ell))$ is $S_{3 \ell-1,1}$ for $q=0$ and $0$ if 
$q \neq 0$.}
\item{ Assume $-9 \le \ell \le -1$. Then  
\[ H^q(\Omega^{-p}(\ell-p)) = 
\begin{cases} 
{\mathbb C}  & \text{for $(p,q)=(\ell,-\ell)$,} \\ 
0            & \text{otherwise.} \end{cases} \] 
Moreover $H^q(\Q(\ell))$ is $S_{-(3\ell+1),-(3\ell+2)}$ if 
$q=2$ and $0$ if $q \neq 2$. }
\end{enumerate} 
Other values of $\ell$ are not difficult to analyse, but we will not 
need them. One can now decompose (\ref{E1pq.EXP}) into irreducible 
pieces using plethysms and the Littlewood-Richardson rule. 

\subsection{Step IV: Analysis of (\ref{neq.SP4}) and $\Hy^\bullet$} 
The term $\acc E_2^{p,q}$ can be nonzero only for $p=-5,-4$,
so $\acc E_3 =\,  \acc E_\infty$. 
By its construction, the differential $\acc d_2$ is a 
composite of two connecting maps 
\[ H^q(\Q(\ell)) \stackrel{\alpha_q}{\lra} H^{q+1}(\O_X(\ell)) 
\stackrel{\beta_q}{\lra} H^{q+2}(\I_X(\ell)) \] 
coming from sequences (\ref{defn.IX}),(\ref{defn.Q}). 
\begin{enumerate}
\item{Assume $\ell \ge 0$, then the groups 
\[ H^{\ge 1}(f_* \O_{\P^2 \times \P^2} \otimes \O_{\P^9}(\ell)), \;\;
   H^{\ge 1}(\O_{\P^9}(\ell)) \] 
vanish. Hence $\alpha_q$ is surjective
for $q=0$ and an isomorphism for $q \ge 1$, whereas $\beta_q$ is 
always an isomorphism. Hence $\acc E_3^{p,q} = 0$ unless 
$(p,q) = (-5,0),(-5,1)$ or $(-4,0)$. Then 
\begin{equation}
 \Hy^{-5} =\, \acc E_3^{-5,0} =\, \acc E_2^{-5,0} = H^0(\I_X(\ell)), 
\label{neq.H5} \end{equation}
and there is an extension 
\begin{equation}
 0 \ra\, \acc E_3^{-5,1} \ra \Hy^{-4} \ra\, \acc E_3^{-4,0} \ra 0. 
\label{neq.H4} \end{equation}
In fact 
$\acc E_3^{-5,1} = \, \acc E_2^{-5,1} = H^1(\I_X(\ell))$ and 
$\acc E_3^{-4,0}$ is the kernel of the map 
$H^0(\Q(\ell)) \ra H^2(\I_X(\ell))$. 

The terms $E_1^{p,q}$ are nonzero only for $q=0$, so $E_2 = E_\infty$ and 
$\Hy^{-5}=E_2^{-5,0},\quad \Hy^{-4}=E_2^{-4,0}$.}

\item{ Now assume $-9 \le \ell \le -1$, then  
the maps $\alpha_q,\beta_q$ are bijective for $0 \le q \le 2$.
The only nonzero $\acc E_2^{-4,q}$ term is at $q=2$, hence 
\begin{equation} 
H^4(\I_X(\ell)) = H^2(\Q(\ell)) = S_{-(3\ell+1),-(3\ell+2)}.
\end{equation}
If  $\ell = -2,-1$, then each $E_1^{p,q}$ term is zero,
hence $\Hy^\bullet$ is identically zero. If $\ell \le -3$, then 
$E_1^{p,q}$ is nonzero exactly when $(p,q)=(\ell,-\ell)$. In that 
case, we have 
\begin{equation} \Hy^0 = S^*_{-(2\ell+3)} \otimes S^*_{-(\ell+3)}, \qquad 
\Hy^j =0 \quad \text{if $j \neq 0$.} 
\end{equation}}
\end{enumerate}

The preparations are over, and the computation proper 
may begin. 
\subsection{Step V} 
The procedure is to use formula (\ref{complex.char}) 
on complexes $\{ {\check E}_1^{p,q} \}_p$ and 
$\{E_1^{p,q}\}_p$ coming from 
the rows of (\ref{ideal.spectralseq}) and (\ref{neq.SP3}). 
As explained earlier, if $\ell >0$, then 
$[E_1^{\bullet,0}]$ can be calculated explicitly 
using plethysms followed by the L-R rule. 

Since $M_{30}$ is a $20$-dimensional summand of $S_3(S_3)$, 
it is necessarily $\{33,30\}$. 
Let $\ell=4$, then 
\begin{equation}
 [\Hy^{-4}] - [\Hy^{-5}] = [E_1^{\bullet,0}] = \sigma_{11,1} - 
(\sigma_{66}+\sigma_{63}+\sigma_{60}+\sigma_{51}+\sigma_{42}
+\sigma_{00}). 
\label{Z1} \end{equation} 
(This is a heavy Maple computation.) 
The only nonzero terms in (\ref{ideal.spectralseq}) are 
${\check E}_1^{0,0}$ and ${\check E}_1^{-1,0}$. Hence 
$H^j(\I_X(4)) = 0 $ for $j \neq 0$ and 
\begin{equation}
 [H^0(\I_X(4))] = [M_{40} \oplus (M_{30} \otimes S_3)] - 
[M_{31}]. 
\label{Z2} \end{equation}
From (\ref{neq.H5}) and (\ref{neq.H4}), 
\begin{equation}
[\Hy^{-5}]  = [H^0(\I_X(4))], \qquad 
[\Hy^{-4}]  = [H^0(\Q(4))] = \sigma_{11,1}. 
\label{Z3} \end{equation}
From $[M_{30}] = \sigma_{33}+ \sigma_{30}$, the value of 
$[M_{30} \otimes S_3]$ is known by the L-R rule. 
From (\ref{Z1})--(\ref{Z3}), we have 
\[ [M_{40}] - [M_{31}] = 
\sigma_{66} - (\sigma_{42}+\sigma_{33}+\sigma_{21}). \] 
Since dim $M_{40} =28$, this forces 
$M_{40} = \{66\}$ and $M_{31} = \{42,33,21\}$.  

Now let $\ell=5$. Then (by an even heavier computation)
\begin{equation} \begin{aligned} 
{} [\Hy^{-4}]-[\Hy^{-5}] =  [E_1^{\bullet,0}] =\, &  
  \sigma_{14,1} - (\sigma_{96}+\sigma_{93}+
\sigma_{90}+ \sigma_{81}+ \sigma_{75} + \\
 & 2\sigma_{72} + \sigma_{63}+ \sigma_{54} 
+\sigma_{51}+ \sigma_{33}+ \sigma_{30}). 
\end{aligned} \label{Y1} \end{equation}
From (\ref{ideal.spectralseq}), we 
deduce that $\I_X(5)$ has no higher cohomology and 
\begin{equation} \begin{aligned} 
{} & [H^0(\I_X(5))] = 
[{\check E}_1^{0,0}] - [{\check E}_1^{-1,0}]
+ [{\check E}_1^{-2,0}]  =  \\ 
& [M_{30} \otimes S_2(S_3)] +  [M_{40} \otimes S_3]
- [M_{31} \otimes S_3] - [M_{41}] + [M_{32}]. 
\end{aligned} \label{Y2} \end{equation}
As before 
\begin{equation}
[\Hy^{-5}]  = [H^0(\I_X(5))], \qquad 
[\Hy^{-4}]  = [H^0(\Q(5))] = \sigma_{14,1}. 
\label{Y3}\end{equation}
Combining (\ref{Y1})--(\ref{Y3}), we have a relation 
\[ [M_{32}] + (\sigma_{75}+\sigma_{54}+\sigma_{33}) = 
   [M_{41}] + (\sigma_{42}+\sigma_{21}+\sigma_{00}),
\]
which implies $M_{32} = \{42,21,00\}$ and 
$M_{41} = \{75,54,33\}$. Evidently we could continue the procedure 
with larger and larger $\ell$, and calculate all the $M_{j,-p}$. 
Alternately, we can twist by a negative 
$\ell$ and work our  way from the other end of the resolution. 

For instance, let $\ell=-1$. The only nontrivial map in 
(\ref{ideal.spectralseq}) is 
\[ {\check E}_1^{-6,9} \lra {\check E}_1^{-5,9}, \quad 
\text{i.e.,} \quad 
M_{46} \otimes S_3^* \lra M_{45}.\]
It has kernel $H^3(\I_X(-1)) = H^1(\Q(-1)) = 0$, and cokernel 
$H^4(\I_X(-1)) = H^2(\Q(-1)) = S_{21}$. Hence 
\[ \sigma_{21} = [M_{45}] - [M_{46} \otimes S_3^*]. 
\] 
Substituting $M_{46} = \{00\}$, we have  $M_{45} = \{ 33, 21\}$. 
Now let $\ell=-2$ to calculate $M_{44}$, etc. 
These are the fruits of our labour: 
\begin{equation} \begin{array}{lll}
M_{30} = \{33,30\}, & M_{31} = \{42,33,21\}, \\ 
M_{32} = \{42,21,00\}, & M_{33} = \{30\}, \\ 
M_{40} = \{66\}, & M_{41} = \{75,54,33\}, \\
M_{42} = \{75,63,54,42,33,21\}, & M_{43} = \{66,63,54,42,42,30,21,00\}, \\
M_{44} = \{54,42,33,30,21\}, & M_{45} = \{33,21\}, \\ 
M_{46} = \{00\}.
\end{array} \label{syzygymodules.neq} \end{equation} 

\begin{Scholium} \rm 
\emph{In principle}, we can work with negative $\ell$ throughout 
and bypass steps II--IV altogether. I see two reasons dictating 
against this. Firstly, the procedure becomes impractical as we 
move rightwards in the resolution. (For instance, in order to reach 
$M_{40}$,  we need to set $\ell=-6$ and then calculate 
$M_{45} \otimes S_5(S_3^*)$ amongst other things.) 
Secondly, (and what is more to the point) I believe that the 
construction used here is of interest not confined to 
this example. In \cite{ego3}, it was used on varieties defined 
by binary forms having roots of specified 
multiplicities. 
\end{Scholium} 
\setcounter{equation}{0}

\section{The Locus $X_\Delta$} 

\begin{equation} \begin{tabular}{c|ccccccccc}
  & 0  & 1   & 2   & 3   & 4   & 5   & 6  & 7  & 8 \\ \hline 
4 & 35 & 119 & 210 & 252 & 210 & 120 & 45 & 10 & 1 \\ 
5 &    &     &  1  &     &     &     &    &    & 
\end{tabular} \label{betti.Delta} \end{equation}

The symmetric group ${\mathfrak S_3}$ acts on $(\P V^*)^{3}$ by 
permuting the factors, and $X_{\Delta}$ is the categorial quotient 
$(\P V^*)^{3}/{\mathfrak S_3}$. Hence for 
$\ell \ge 0$, the group $H^0(\O_X(\ell))$ is the 
${\mathfrak S_3}$--invariant part of 
\[ (H^0(\P V^*,\O_{\P}(\ell)))^{\otimes 3} = (S_\ell)^{\otimes 3},\] 
which is $S_3(S_\ell)$. By the same argument, 
$H^j(\O_X(\ell)) = 0$ for $\ell \in {\mathbf Z},0 < j < 6$; and 
$H^6(\O_X(\ell)) = S_3(S_{-(\ell+3)}^*)$ for $\ell \le -3$. 
Now we have identifications 
\[ \begin{array}{cll}
0 \lra  H^0(\I_X(\ell)) \lra & H^0(\P^9,\O_{\P}(\ell)) & 
\lra H^0(X,\O_X(\ell)) \\
 & \quad \;\;  \big| \big| & \qquad \quad \big| \big| \\
 & S_\ell(S_3) & \stackrel{\gamma}{\lra}  S_3(S_\ell)
\end{array} \]

If $\ell \ge 3$, then in (\ref{ideal.spectralseq}) we have 
${\check E}_1^{p,q} =0$ for $q > 0$, so $H^1(\I_X(\ell))=0$. 
Thus $\gamma$ is surjective for $\ell \ge 3$. (This is a rather 
special instance of Foulkes's conjecture--see 
e.g. \cite[p.141]{MacDonald}.) 
Now formula (\ref{E1.equiv}) holds, and the calculation 
proceeds as in that case. This gets tedious for high values 
of $\ell$, but one can bypass it by calculating from 
the other end of the resolution. 

For example, let $\ell=1$. The only nontrivial map in 
(\ref{ideal.spectralseq}) is 
\[ 
{\check E}_1^{-8,9} \stackrel{d_1}{\lra} {\check E}_1^{-7,9}, 
\quad \text{i.e.,} \quad 
 M_{48} \otimes H^9(\O_{\P}(-11)) \lra M_{47}. \]
Since $H^j(\P^9, \I_X(1)) = 0$ for $j=1,2$, this map is an isomorphism.   
Hence $M_{47} = \{33\}$. Now let $\ell=0$ and calculate $M_{46}$, etc. 
The outcome is 
\begin{equation} \begin{array}{l}
M_{40} = \{51\},  \, M_{41} = \{63,60,42\}, 
\, M_{42} = \{72,66,63,42,30\}, \\ 
M_{43} = \{75,72,54,51,33,30\}, \, M_{44} = \{75,63,60,42,33\}, \\
M_{45} = \{66,63,42,00\}, \, M_{46} = \{54,30\}, \, M_{47} = \{33\}, \\
M_{48} = \{00\}, \, M_{52} = \{00\}. 
\end{array} \label{syzygymodules.Delta} 
\end{equation}

In \cite[\S 226]{Salmon1}, Salmon gives a system of forty-five 
equations each of degree $4$ in the $\abar$, which is satisfied 
iff the cubic $F$ breaks up into three lines. Thus the system 
generates an ideal $J$ such that $\sqrt{J} = I_X$. His construction 
is natural with respect to the $SL_3$-action, hence $(J)_4 \subseteq M_{40}$ 
is a subrepresentation. Since $M_{40}$ is irreducible, $J = I_X$. It 
follows that there are ten dependencies in his system, something which is not 
a priori evident. 

\begin{Problem} \rm 
\begin{enumerate}
\item{The resolution shows some measure of duality, viz.,
$ M_{42}=M_{44}^*, \, M_{41} \oplus \complex = M_{45}^*, \,
M_{43} = M_{43}^*$. 
There ought to be a theoretical explanation.} 
\item{It follows from simple degree considerations  that $M_{40}$ must 
be identical to the set of Brill's equations (\cite[p.139 ff]{GKZ}). It 
would be worthwhile to check this by direct calculation. }
\end{enumerate}\end{Problem} 
\setcounter{equation}{0}

The next two computations are less conceptual that the preceding 
ones, in that they rely upon incidental features of the 
resolution. 

\section{The Locus $X_\emptyset$}

\begin{equation} \begin{tabular}{c|cccc}
  & 0  & 1  & 2  & 3  \\ \hline 
8 & 35 & 70 & 45 & 8  \\
9 &    &    &    & 1 
\end{tabular} \label{betti.emptyset} \end{equation} 

We decompose  $S_8(S_3)$ and look for a $35$-dimensional 
subrepresentation. The only possibilities 
for $M_{80}$ are found to be $\{51\}$ or $\{54\}$. 
In fact, it turns out that $M_{80} = \{54\}$. Let us grant this 
for the moment and postpone the justification to 
\S \ref{concom.emptyset}. 

Substitute $\ell =9$ in (\ref{ideal.spectralseq}). The only 
nontrivial map ${\check E}_1^{-1,0} \stackrel{d_1}{\lra} 
{\check E}_1^{0,0}$, i.e., 
$M_{81} \ra M_{80} \otimes S_3$, is necessarily injective. 
On dimensional grounds $M_{81} = \{54,42,21\}$. Now 
$H^6(\O_X)=0$ (this was verified in Macaulay-2), hence 
$H^7(\I_X)=0$. But from (\ref{ideal.spectralseq}), this 
group is the cokernel of the map 
\[ (M_{83} \otimes S_3^*) \, \oplus S_2(S^*_3) \lra 
   M_{82}. \]
Since $S_2(S_{33}) = \{66,42\}$, it cannot by itself 
contribute a $45$-dimensional summand. But then 
$M_{83} = \{21\}$ is forced. Now the only possibilities 
for $M_{82}$ are $\{54,33\},\{42,33,21\}$. Since the resolution is 
minimal, the map 
$M_{93} \otimes R(-12) \lra M_{82} \otimes R(-10)$ is 
nonzero. Hence $M_{82} \otimes R_2 = M_{82} \otimes S_2(S_3)$ 
must contain a trivial summand. This rules out 
$\{54,33\}$, so $M_{82} = \{42,33,21 \}$. In 
sum, 
\begin{equation} \begin{array}{lll}
M_{80} = \{54\}, & M_{81} = \{54,42,21\}, & M_{82} = \{42,33,21\}, \\
M_{83} = \{21\}, & M_{93} = \{00\}. 
\end{array} 
\label{syzygymodules.emptyset} \end{equation}

In \cite[\S 240]{Salmon1}, Salmon gives a system of equations each of degree 
$8$ in the $\abar$, which is satisfied iff the cubic $F$ is factorisable. 
Using an argument similar to the case $X_\Delta$, it follows that his system 
generates $I_{X_\emptyset}$.

\begin{Problem} \rm 
Define $\Q = \coker \, (\O_{X_\emptyset} \ra 
f_* \O_{\P S_2^* \times \P V^*})$ as in the case of $X_{\neq}$. 
It would be of interest to calculate it explicitly along the lines 
of lemma \ref{lemma.calc.Q}, or even to calculate the 
groups $H^\bullet(\Q(\ell))$. Once this is done, there is no 
difficulty in adapting steps II--V to $X_{\emptyset}$. 
\end{Problem} 
\setcounter{equation}{0}

\section{The Locus $X_{\tact}$} 

\begin{equation} \begin{tabular}{c|ccc} 
  & 0  & 1  & 2  \\ \hline 
4 & 1          \\ 
5 & 10 & 8     \\ 
6 &    & 10 & 8 
\end{tabular} \label{betti.tact} \end{equation}  

Now $M_{40}=\{00\}$ and 
$H^0(\I_X(5)) = (M_{40} \otimes S_3) \oplus M_{50}$ is a 
$20$--dimensional module inside $H^0(\O_{\P}(5)) = S_5(S_3)$. 
Decomposing the latter, on dimensional grounds the module 
can only be $\{33,30\}$. So $M_{50} = \{33\}$. 

Now $M_{51}$ is a summand of $(M_{40} \otimes R_2) \oplus 
(M_{50} \otimes R_1)$, so it must be $\{21\}$. Since 
$M_{61}$ is a summand of 
$(M_{40} \otimes R_3)\oplus (M_{50} \otimes R_2)$, it is 
either $\{33\}$ or $\{30\}$, we do not yet know which. 
Decompose $(M_{61} \otimes R_1) \oplus  
(M_{51} \otimes R_2)$ for each of the possibilities, and one 
sees that $M_{62}$ is necessarily $\{21\}$. 

It remains to determine $M_{61}$. Since $X$ is $6$-dimensional, 
$H^8(\I_X(-3))=0$. Substituting $\ell=-3$ in 
(\ref{ideal.spectralseq}), this group appears as the cokernel of 
\[M_{62} \otimes H^9(\P^9,\O_{\P}(-11)) \lra M_{61}.\]
The only $10$-dimensional summand of $\{21\} \otimes \{33\}$ is 
$\{33\}$, so $M_{61}=\{33\}$. Thus 
\begin{equation}
M_{40} = \{00\}, M_{50} = \{33\},  M_{51} = \{21\}, 
M_{61} = \{33\},  M_{62} = \{21\}. 
\label{syzygymodules.tact} \end{equation}
\begin{Problem} \rm 
Without recourse to a machine computation, show that 
$X_{\tact}$ is arithmetically Cohen--Macaulay. (This will 
probably involve the theory of complete conics.) 
\end{Problem} 
\setcounter{equation}{0}
\section{The Locus $X_Y$} \label{section.XY} 
\begin{equation}
\begin{tabular}{c|cccc} 
  & 0  & 1  & 2  & 3 \\ \hline 
3 & 20 & 45 & 36 & 10 
\end{tabular} 
\label{betti.Y} \end{equation} 

This is a rank variety in the sense of Porras \cite{Porras}, 
and its resolution is deduced there (Prop. 4.2.3). 
We summarise the solution.  

A subspace $W$ of $V^*$ gives an inclusion 
$\P (S_3 W) \subseteq  \P (S_3 V^*)$. Then 
\[ X_Y = \{ F: F \in \P(S_3 W)\; \text{for some \emph{proper}
subspace $W$}\}.\]
(This is so, because the $L_i$ define concurrent lines iff they fail to span 
$V^*$.)
The multiplication $S_2 \otimes S_1 \lra S_3$ gives a map 
of vector bundles 
\[ \alpha: \underbrace{S_2 \otimes \O_{\P}(-1)}_{\bb}
\lra V^* \otimes \O_{\P} \] 
on $\P S_3^*$. 
Then $X$ coincides (as a scheme) with the 
degeneracy locus $\{ \rank\, \alpha \le 2 \}$, and 
the Eagon--Northcott complex of $\alpha$ resolves 
its structure sheaf. 
\[
0 \ra \wedge^6 \bb \otimes S_3 \ra 
\wedge^5 \bb \otimes S_2 \ra \wedge^4 \bb \otimes S_1 \ra 
\wedge^3 \bb \ra \O_{\P^9} \ra \O_{X_Y} \ra 0.
\]
Now $\wedge^j \bb = \wedge^j S_2 \otimes \O_{\P}(-j)$, 
and after decomposing these further
\begin{equation} 
M_{30} = \{33,30\}, 
M_{31} = \{42,33,21\}, 
M_{32} = \{42,21,00\}, 
M_{33} = \{30\}. 
\label{syzygymodules.Y} \end{equation}

\section{Ideal Generators and Concomitants} 
Let $X$ be any of the loci above. For $\ell \ge 1$, we have an 
exact sequence 
\[  (I_X)_{\ell-1} \otimes S_3 \lra S_{\ell}(S_3) 
\stackrel{g_\ell}{\lra} U_{\ell} \lra 0. 
\] 
The module $M_{\ell,0}$ consists of degree $\ell$ primitive 
generators for $I_X$, as such it is a submodule of 
$U_\ell$. Let $S_{m+n,n} \subseteq M_{\ell,0}$ be 
a direct summand. By choosing an equivariant splitting of $g_\ell$, 
we may write $S_{m+n,n} \subseteq S_{\ell}(S_3)$. Then 
by \S \ref{section.concomitants}, to specify this inclusion 
is to specify a concomitant $\Phi_{\ell,m,n}$. 

It turns out that in \emph{every} case under consideration, there 
is only one copy of $S_{m+n,n}$ inside $S_\ell(S_3)$, hence 
the inclusion is independent of the choice of the 
splitting.\footnote{There is no difficulty in adapting to the 
case when this is not so.} Thus 
it is enough to produce some symbolic 
expression of type $(\ell,m,n)$ (by trial and error), and to verify 
that it does not vanish identically after substitution. Then 
it must be the $\Phi_{\ell,m,n}$ that we are looking for. Besides 
doing this, in each case I have verified by direct computation 
(in Macaulay-2) that $\Phi$ vanishes on $X$. Although 
logically superfluous, this is a very useful check against 
the propagation of error. 

\subsection{Case $X_{\equiv}$} \label{concomitants.equiv}
Corresponding to $M_{20} = \{42\}$, we need a concomitant 
of type $(2,2,2)$. The only legal symbolic expression is 
\[ 
\Phi_{222} = (\alpha \, \beta \, u)^2\, \alpha_x \, \beta_x. \]

To recapitulate, let $T_1,\dots,T_{27}$ be the semistandard 
tableaux on $(4,2)$ (see \S \ref{basis.Smn}) and write 
$ \Phi_{222} = 
\sum\limits_{i=1}^{27} f_i \otimes X_{T_i}$. 
(A priori, after expansion $\Phi$ will contain nonstandard 
monomials such as $x_1^2u_1u_3$. These can be rewritten as linear 
combinations of the $X_T$ using the straightening 
law--see \cite[\S 8.1]{Fulton}.)
Then $I_X = (f_1,\dots,f_{27})$. This can be rephrased as a criterion: 

{\it A ternary cubic $F$ can be written as the cube of a 
linear form iff the concomitant 
$(\alpha \, \beta \, u)^2\, \alpha_x \, \beta_x$ vanishes on 
$F$.} There is a similar claim for each $X$. 

\subsection{Case $X_{\neq}$}

\[ M_{30} \longmapsto 
\Phi_{303} = \brac{\alpha}{\beta}{\gamma} 
\brac{\alpha}{\beta}{u} \brac{\alpha}{\gamma}{u} 
\brac{\beta}{\gamma}{u},\; 
\Phi_{330} = {\brac{\alpha}{\beta}{\gamma}}^2 
             \alpha_x \,\beta_x \, \gamma_x.  
\]
These are respectively the Cayleyan and the Hessian of $F$ 
(see \cite[\S 218,219]{Salmon1}).  
\[ M_{40} \longmapsto 
\Phi_{406} = {\brac{\alpha}{\beta}{u}}^2 {\brac{\gamma}{\delta}{u}}^2
              \brac{\alpha}{\delta}{u} \brac{\beta}{\gamma}{u}. 
  \] 

\subsection{Case $X_\Delta$}
\[ M_{40} \longmapsto \Phi_{441} = 
\brac{\alpha}{\beta}{\gamma}\,\brac{\alpha}{\gamma}{\delta}\,
\brac{\alpha}{\beta}{u}\, \beta_x\,\gamma_x\,\delta_x^{\,2}.  
\] 

\subsection{Case $X_\emptyset$} 
\label{concom.emptyset} 
We are yet to prove that $M_{80}$ is $\{54\}$ and not 
$\{51\}$. 
There is a well known invariant of ternary cubics 
in degree 4, namely the Aronhold invariant (see \cite[p.~167]{Sturmfels}) 
\[ \Phi_{400} = 
\brac{\alpha}{\beta}{\gamma}\brac{\alpha}{\beta}{\delta}
\brac{\alpha}{\gamma}{\delta}\brac{\beta}{\gamma}{\delta}. 
\] 
Decomposing $S_8(S_3)$, one sees that it houses exactly one copy each of 
$S_{51},S_{54}$. Hence there is one 
concomitant each of type $(8,4,1),(8,1,4)$. Now $\Phi_{841}$ is necessarily 
equal to the product $\Phi_{400}\Phi_{441}$. Since $I_{X_{\emptyset}}$ has 
no generators in degree $4$, none of the factors can vanish on 
$X_{\emptyset}$. Hence $M_{80} \neq \{51\}$. Thus 
\[ 
M_{80} \longmapsto 
\Phi_{814} = 
\alpha_x \brac{\alpha}{\beta}{\gamma} 
\brac{\alpha}{\beta}{\delta} \brac{\beta}{\gamma}{\epsilon} 
\brac{\gamma}{\zeta}{u} \brac{\delta}{\epsilon}{u} 
\brac{\delta}{\eta}{u} \brac{\epsilon}{\theta}{u} 
{\brac{\zeta}{\eta}{\theta}}^2. 
\] 
It is quite easy to concoct symbolic expressions of type 
$(8,1,4)$, but most will vanish identically after expansion. 
It cost the author a week's labour to find one which does not. 

\subsection{Case $X_{\tact}$}
\[ M_{40} \longmapsto \Phi_{400}\;\, (\text{as written above}). \] 
\[ M_{50} \longmapsto \Phi_{503} = 
\brac{\alpha}{\beta}{\gamma}\,
\brac{\alpha}{\beta}{\delta}\,\brac{\beta}{\gamma}{\epsilon}\,
\brac{\alpha}{\gamma}{u} \,{\brac{\delta}{\epsilon}{u}}^2. \] 

\subsection{Case $X_Y$}
\[ M_{30} \longmapsto 
   \Phi_{303},\Phi_{330} \;\, (\text{as written above}). 
\] 
In fact, the Hessian vanishes identically iff the Cayleyan 
does. (This is an easy deduction from the formulae in 
\cite{Cayley1}.) Hence either one defines the locus set-theoretically. 
Also see \cite[p.~234]{Olver} for a discussion of the 
Gordan--N{\"o}ther theorem.

\begin{Scholium} \rm The symbolic method can be used for describing the higher 
syzygies as well. We will give an illustration of this idea but will not attempt any 
systematic treatment. Consider the resolution of $X_{\equiv}$, in particular the 
submodule $S_{54} \subseteq M_{21}$. Keep the notation of 
\S \ref{concomitants.equiv}. We have the differential 
\[ S_{54} \otimes \O_{\P}(-3) \stackrel{\partial^1}{\lra}
 S_{42} \otimes \O_{\P}(-2), \] 
or equivalently, a map 
\[ \O_{\P} \lra S_{42} \otimes S_{51} \otimes\O_{\P}(1), \] 
or equivalently, an element 
\[ \Psi \in H^0(\P^9,S_{42} \otimes S_{51} \otimes \O_{\P}(1)) = 
S_{42} \otimes S_{51} \otimes S_3. \] 
We would like to represent $\Psi$ symbolically. Let $\ybar = \{y_1,y_2,y_3\}$ be 
`copies' of the variables $\{x_1,x_2,x_3\}$, and ditto for 
$\vbar = \{v_1,v_2,v_3\}$ answering to $\{u_1,u_2,u_3\}$. 
Let $\TT_1, \dots, \TT_{35}$ be the semistandard tableaux on numbers $1,2,3$ on 
the Young diagram of $(5,1)$. Exactly as in \S \ref{basis.Smn}, $Y_{\TT_j}$ will 
denote a monomial in the $\ybar,\vbar$ corresponding to $\TT_j$.
Write  $\alpha_y = \alpha_1 y_1 + \alpha_2 y_2 + \alpha_3 y_3, \, 
v_x = v_1 x_1 + v_2 x_2 + v_3 x_3$, and 
\[ 
\brac{\alpha}{u}{v} = \left| \begin{array}{ccc} 
\alpha_1 & \alpha_2 & \alpha_3 \\ 
u_1 & u_2 & u_3 \\ 
v_1 & v_2 & v_3 \end{array} \right|. 
\] 
A calculation (which we suppress) shows that 
$\Psi = \brac{\alpha}{u}{v}^2 \alpha_y v_x^2$. That is to say, after expanding 
the product and substituting, we can write 
\[ \Psi = \sum\limits_{j=1}^{35} \; 
(\, \sum_{i=1}^{27} h_{ij} \otimes X_{T_i})
\otimes Y_{\TT_j},
\] 
where $h_{ij}$ are linear forms in the $\abar$. Let $\Sigma_j$ be the 
sum in the parentheses multiplying $Y_{\TT_j}$. Then $\Sigma_j$ represents a 
linear syzygy between the ideal generators $f_i$, namely 
$\sum\limits_{i=1}^{27} h_{ij} f_i =0$. This accounts for the $35$ syzygies 
corresponding to the summand $S_{54} \subseteq M_{21}$, and of course 
other summands are given by other symbolic expressions. In fact 
\[ S_{51} \longmapsto \alpha_x \alpha_y^2 v_x u_y^2, \quad 
   S_{42} \longmapsto {\brac{\alpha}{u}{v}} \alpha_x \alpha_y v_x u_y, \quad 
   S_{21} \longmapsto {\brac{\alpha}{u}{v}} \alpha_x^2 u_y.   
\] 
\end{Scholium} 
\smallskip 

The complete system of concomitants for a ternary cubic 
is given by Cayley \cite{Cayley1}. I enclose a table collating our notation 
with his. 
\[ \begin{array}{llllll}
\Phi_{222} = \Theta & (\text{no.$4$}), & 
\Phi_{303} =  P     & (\text{no.$5$}), & 
\Phi_{330} =  H     & (\text{no.$12$}), \\
\Phi_{406} =  F     & (\text{no.$17$}), &
\Phi_{441} = J      & (\text{no.$19$}), & 
\Phi_{400} = S      & (\text{no.$1$}), \\
\Phi_{814} = \overline{J} & (\text{no.$27$}), & 
\Phi_{503} = Q      & (\text{no.$16$}). & 
\end{array} \] 

\begin{Scholium} \rm 
The invariant $S$ is ubiquitous in the geometry of plane 
cubics, it is essentially the same as the Eisenstein series $g_2$. There is 
a sextic invariant (also called $T$)
\[ \Phi_{600} = 
\brac{\alpha}{\beta}{\gamma} \brac{\alpha}{\beta}{\delta} 
\brac{\beta}{\gamma}{\epsilon} 
\brac{\alpha}{\gamma}{\zeta} \brac{\delta}{\epsilon}{\zeta}^2, \] 
which is $g_3$ in a different guise. The discriminant $\Delta$ of $F$ is 
a linear combination of $S^3$ and $T^2$ (the coefficients depend on one's 
normalisation). Now $\Delta =0$ is the closure of the locus of nodal cubics. 
The locus of cuspidal cubics is the complete intersection 
$S=T=0$. (See \cite[\S 224]{Salmon1}.)
\end{Scholium}

Several geometric theorems about the Cayleyan can be found in \cite{Cayley2} (where 
it is called the Pippian). The invariant $S$ vanishes on a 
smooth cubic curve iff the latter is projectively isomorphic to the Fermat 
cubic $x_1^3 + x_2^3 + x_3^3 =0$. Such cubics are called 
anharmonic. See \cite{DolgachevKanev} for some beautiful results in this 
vein. 

\bibliographystyle{plain}
\bibliography{references.cubics}

\bigskip 

\parbox{13cm}{\small 
Jaydeep V. Chipalkatti \\
Department of Mathematics and Statistics, \\
416 Jeffery Hall, Queen's University, \\ 
Kingston, ON K7L 3N6, Canada. \\ 
email: {\tt jaydeep@mast.queensu.ca}}

\end{document}